\documentclass[11pt]{article}
\usepackage{amsmath}
\usepackage{amssymb}
\usepackage{graphicx}
\DeclareMathOperator{\sgn}{sgn}
\DeclareMathOperator{\diag}{diag}
\def \v{\mathbf{v}}
\def \u{\mathbf{u}}
\def \w{\mathbf{w}}
\def \x{\mathbf{x}}
\def \y{\mathbf{y}}
\def \i{\mathbf{i}}
\def \A{\mathcal{A}}
\def \L{\mathcal{L}}
\def \Q{\mathcal{Q}}
\def \lamax{\lambda_{\max}}

\def \e{\epsilon}
\def \la{\lambda}
\def \T{\mathcal{T}}
\def \Go{G^{\large \mathtt  o}}
\def \E{\mathcal{E}}
\def \D{\mathcal{D}}
\newtheorem{theorem}{\scshape \mdseries  Theorem}[section]
\newtheorem{defi}[theorem]{\scshape \mdseries  Definition}
\newtheorem{lemma}[theorem]{\scshape \mdseries  Lemma}
\newtheorem{coro}[theorem]{\scshape \mdseries  Corollary}
\newtheorem{conj}[theorem]{\scshape \mdseries  Conjecture}

\begin{document}

\title{\sf The largest $H$-eigenvalue and spectral radius of Laplacian tensor of non-odd-bipartite generalized power hypergraphs\thanks{
Supported by National Natural Science Foundation of China (11371028),
Natural Science Research Foundation of Anhui Provincial Department of Education (KJ2015A322),
Scientific Research Fund for Fostering Distinguished Young Scholars of Anhui
University (KJJQ1001), Project of Academic Innovation Team of Anhui University (KJTD001B),
Open Project of School of Mathematical Sciences of Anhui University (ADSY201501). }
}
\author{Yi-Zheng Fan$^{1,}$\thanks{Corresponding author.
 E-mail addresses: fanyz@ahu.edu.cn(Y.-Z. Fan), muradulislam@foxmail.com (M. Khan), tansusan1@ahjzu.edu.cn (Y.-Y. Tan)},~ Murad-ul-Islam Khan$^1$,  Ying-Ying Tan$^{1,2}$\\
       {\small  \it $1$. School of Mathematical Sciences, Anhui University, Hefei 230601, P. R. China}\\
    {\small  \it $2$. Department of Mathematics and Physics, Anhui Jianzhu University, Hefei 230601, P. R. China}
 }
\date{}
\maketitle

{\small
\noindent
\textbf{Abstract}:
Let $G$ be a simple graph or hypergraph, and let $\A(G),\L(G),\Q(G)$ be the adjacency, Laplacian and signless Laplacian tensors of $G$ respectively.
The largest $H$-eigenvalues (respectively, the spectral radii) of $\L(G),\Q(G)$ are denoted respectively by $\lamax^\L(G), \lamax^\Q(G)$
(respectively, $\rho^\L(G), \rho^\Q(G)$).
It is known that for a connected non-bipartite simple graph $G$, $\lamax^\L(G)=\rho^\L(G) < \rho^\Q(G)$.
But this does not hold for non-odd-bipartite hypergraphs.
We will investigate this problem by considering a class of generalized power hypergraphs $G^{k,\frac{k}{2}}$, which are constructed from simple connected graphs $G$ by blowing up each vertex of $G$ into a $\frac{k}{2}$-set and preserving the adjacency of vertices.

Suppose that $G$ is non-bipartite, or equivalently $G^{k,\frac{k}{2}}$ is non-odd-bipartite.
We get the following spectral properties: (1) $\rho^\L(G^{k,{k \over 2}}) = \rho^\Q(G^{k,{k \over 2}})$ if and only if $k$ is a multiple of $4$; in this case $\lamax^\L(G^{k,\frac{k}{2}})<\rho^\L(G^{k,\frac{k}{2}})$.
(2) If $k\equiv 2 (\!\!\!\mod 4)$, then for sufficiently large $k$, $\lamax^\L(G^{k,\frac{k}{2}})<\rho^\L(G^{k,\frac{k}{2}})$.
Motivated by the study of hypergraphs $G^{k,\frac{k}{2}}$, for a connected non-odd-bipartite hypergraph $G$, we give a characterization of $\L(G)$ and $\Q(G)$ having the same spectra or the spectrum of $\A(G)$ being symmetric with respect to the origin, that is, $\L(G)$ and $\Q(G)$, or $\A(G)$ and $-\A(G)$ are similar via a complex (necessarily non-real) diagonal matrix with modular-$1$ diagonal entries.
 So we give an answer to a question raised by Shao et al., that is, for a non-odd-bipartite hypergraph $G$, that $\L(G)$ and $\Q(G)$ have the same spectra can not imply they have the same $H$-spectra.


}

\section{Introduction}
A {\it hypergraph} $G=(V(G),E(G))$ consists of a set of vertices say $V(G)=\{v_1,v_2,\ldots,v_n\}$ and a set of edges say $E(G)=\{e_{1},e_2,\ldots,e_{m}\}$
  where $e_{j}\subseteq V(G)$.
If $|e_{j}|=k$ for each $j=1,2,\ldots,m$, then $G$ is called a {\it $k$-uniform} hypergraph.
In particular, the $2$-uniform hypergraphs are exactly the classical simple graphs.
For a $k$-uniform hypergraph $G$, if we add to $G$ some edges with cardinality less than $k$, the resulting hypergraph denoted by $\Go$ is one with loops;
and those added edges are called the {\it loops} of $\Go$.
The {\it degree} $d_v(G)$ or simply $d_v$ of a vertex $v \in V(G)$ is defined as $d_v(G)=|\{e_{j}:v\in e_{j}\in E(G)\}|$.
So, a loop contributes $1$ to the degree of the vertex to which it is attached.

An  even uniform hypergraph $G$ is called {\it odd-bipartite} if $V(G)$ has a bipartition $V(G)=V_{1}\cup V_{2} $
   such that each edge has an odd number of vertices in both $V_{1}$ and $V_{2}$.
Hu, Qi and Shao \cite{HQS} introduced the {\it cored hypergraphs} and the {\it power hypergraphs},
   where the cored hypergraph is one such that each edge contains at least one vertex of degree $1$,
   and the $k$-th power of a simple graph $G$, denoted by $G^k$, is obtained from $G$ by replacing each edge (a $2$-set) with a $k$-set by adding $(k-2)$ new vertices.
These two kinds of hypergraphs are both odd-bipartite.

Recently spectral hypergraph theory  has emerged as an important field in algebraic graph theory.
Let $G$ be a $k$-uniform hypergraph on $n$ vertices $v_1,v_2,\ldots,v_n$.
The {\it adjacency tensor} $\A(G)$ of $G$ is defined as $\mathcal{A}(G)=(a_{i_{1}i_{2}\ldots i_{k}})$, a $k$th order $n$-dimensional symmetric tensor, where
 $  a_{i_{1}i_{2}\ldots i_{k}}=\frac{1}{(k-1)!}$ if $\{v_{i_{1}},v_{i_{2}},\ldots,v_{i_{k}}\} \in E(G)$ and $  a_{i_{1}i_{2}\ldots i_{k}}=0$ otherwise.
 Let $\mathcal{D}(G)$ be a $k$th order $n$-dimensional diagonal tensor,
    where $d_{i\ldots i}=d_{v_i}$ for all $i \in [n]:=\{1,2,\ldots,n\}$.
Then $\L(G)=\mathcal{D}(G)-\A(G)$ is called the {\it Laplacian tensor} of $G$,
and $\Q(G)=\mathcal{D}(G)+\A(G)$ is called the {\it signless Laplacian tensor} of $G$.

For a hypergraph $\Go$ with loops, the adjacency tensor of $\Go$ is defined as the same as that of $G$, i.e. $\A(\Go)=\A(G)$.
The Laplacian tensor and the signless Laplacian tensor are defined by $\L(\Go)=\mathcal{D}(\Go)-\A(G)$ and $\Q(\Go)=\mathcal{D}(\Go)+\A(G)$, respectively.
So, even if $\Go$ is not uniform, the adjacency, Laplacian and signless Laplacian tensor of $\Go$ are all $k$th order $n$-dimensional tensors.

In general,
  a real {\it tensor} (also called {\it hypermatrix}) $\mathcal{T}=(t_{i_{1}\ldots i_{k}})$ of order $k$ and dimension $n$ refers to a
  multidimensional array with entries $t_{i_{1}\ldots i_{k}}$ such that $t_{i_{1}\ldots i_{k}}\in \mathbb{R}$ for all $i_{j}\in [n]$ and $j\in [k]$.
 The tensor $\mathcal{T}$ is called \textit{symmetric} if its entries are invariant under any permutation of their indices.
 A {\it subtensor} of $\mathcal{T}$ is a multidimensional array with entries $t_{i_{1}\ldots i_{k}}$ such that $i_j \in S_j \subseteq [n]$ for some $S_j$'s and $j \in [k]$,
denoted by $\mathcal{T}[S_1|S_2|\cdots|S_k]$.
If $S_1=S_2=\cdots=S_k=:S$, then we simply write $\mathcal{T}[S_1|S_2|\cdots|S_k]$ as $\mathcal{T}[S]$, which is called the {\it principal subtensor} of $\mathcal{T}$.
If $k=2$, then $\mathcal{T}[S]$ is exactly the principal submatrix of $\mathcal{T}$;
and if $k=1$, then $\mathcal{T}[S]$ is the subvector of $\mathcal{T}$.

 Given a vector $x\in \mathbb{R}^{n}$, $\mathcal{T}x^{k}$ is a real number, and $\mathcal{T}x^{k-1}$ is an $n$-dimensional vector, which are defined as follows:
   $$\mathcal{T}x^{k}=\sum_{i_1,i_{2},\ldots,i_{k}\in [n]}t_{i_1i_{2}\ldots i_{k}}x_{i_1}x_{i_{2}}\cdots x_{i_k},~
   (\mathcal{T}x^{k-1})_i=\sum_{i_{2},\ldots,i_{k}\in [n]}t_{ii_{2}\ldots i_{k}}x_{i_{2}}\cdots x_{i_k} \mbox{~for~} i \in [n].$$
 Let $\mathcal{I}$ be the {\it identity tensor} of order $k$ and dimension $n$, that is, $i_{i_{1}i_2 \ldots i_{k}}=1$ if and only if
   $i_{1}=i_2=\cdots=i_{k} \in [n]$ and $i_{i_{1}i_2 \ldots i_{k}}=0$ otherwise.

\begin{defi}{\em \cite{Qi2}} Let $\mathcal{T}$ be a $k$th order $n$-dimensional real tensor.
For some $\lambda \in \mathbb{C}$, if the polynomial system $(\lambda \mathcal{I}-\mathcal{T})x^{k-1}=0$, or equivalently $\mathcal{T}x^{k-1}=\lambda x^{[k-1]}$, has a solution $x\in \mathbb{C}^{n}\backslash \{0\}$,
then $\lambda $ is called an eigenvalue of $\mathcal{T}$ and $x$ is an eigenvector of $\mathcal{T}$ associated with $\lambda$,
where $x^{[k-1]}:=(x_1^{k-1}, x_2^{k-1},\ldots,x_n^{k-1})$.
\end{defi}

If $x$ is a real eigenvector of $\mathcal{T}$, surely the corresponding eigenvalue $\lambda$ is real.
In this case, $x$ is called an {\it $H$-eigenvector} and $\lambda$ is called an {\it $H$-eigenvalue}.
The {\it spectral radius of $\T$} is defined as
$\rho(\T)=\max\{|\lambda|: \lambda \mbox{ is an eigenvalue of } \T \}.$
Denote respectively the largest $H$-eigenvalues (respectively, the spectral radii) of $\A(G),\L(G),\Q(G)$ by $\lamax^\A(G),\lamax^\L(G), \lamax^\Q(G)$
 (respectively, $\rho^\A(G),\rho^\L(G), \rho^\Q(G)$).
By Perron-Frobenius theorem of nonnegative tensors (see \cite{CPZ, FGH, YY}), $\lamax^\A(G)=\rho^\A(G)$, $\lamax^\Q(G)=\rho^\Q(G)$.
But this does not hold for the Laplacian tensors in general.

Qi \cite{Qi} showed that $\rho^\L(G) \le \rho^\Q(G)$, and posed a question of identifying the conditions under which the equality holds.
So $$\lamax^\L(G) \le \rho^\L(G) \le \rho^\Q(G)=\lamax^\Q(G).\eqno(1.1)$$
Hu et al. \cite{HQX} proved the following result.

\begin{theorem}{\em \cite{HQX}}\label{Hu}
Let $G$ be a connected $k$-uniform hypergraph.
Then $\lamax^\L(G)=\lamax^\Q(G)$ if and only if $k$ is even and $G$ is odd-bipartite.
\end{theorem}


Denote by $Spec(\A(G))$, $Spec(\L(G))$ and $Spec(\Q(G))$ the spectra of $\A(G)$, $\L(G)$ and $\Q(G)$ respectively, and by
$Hspec(\L(G))$, $Hspec(\L(G))$ and $Hspec(\Q(G))$ the sets of distinct $H$-eigenvalues of $\A(G)$, $\L(G)$ and $\Q(G)$ respectively.
Shao et al. \cite{SSW} gave some characterizations on these different types of spectra.

\begin{theorem}\label{shao0} {\em  \cite{SSW}}
Let $G$ be a connected $k$-uniform hypergraph.
Then $\rho^\L(G) = \rho^\Q(G)$ if and only if $Spec(\L(G))=Spec(\Q(G))$.
\end{theorem}

\begin{theorem}\label{shao} {\em  \cite{SSW}}
Let $G$ be a connected $k$-uniform hypergraph.
Then the following conditions are equivalent.

$(1)$ $k$ is even and $G$ is odd-bipartite.

$(2)$ $Spec(\L(G))=Spec(\Q(G))$ and $Hspec(\L(G))=Hspec(\Q(G))$.

$(3)$ $Hspec(\L(G))=Hspec(\Q(G))$.

$(4)$ $Spec(\A(G))=-Spec(\A(G))$ and $Hspec(\A(G))=-Hspec(\A(G))$, i.e. both $Spec(\A(G))$ and $Hspec(\A(G))$ are symmetric with respect to the origin.

$(5)$ $Hspec(\A(G))=-Hspec(\A(G))$.

\end{theorem}

Suppose that $k$ is even and $G$ is connected.
If $G$ is odd-bipartite, then $\lamax^\L(G)=\lamax^\Q(G)$, which implies that $\lamax^\L(G)=\rho^\L(G)$.
Suppose that $G$ is non-odd-bipartite. Then $\lamax^\L(G)<\lamax^\Q(G)$.
From the inequalities in (1.1), we want to know under which condition
 $\rho^\L(G) = \rho^\Q(G)$ or $\lamax^\L(G)=\rho^\L(G)$.
 If $\rho^\L(G) = \rho^\Q(G)$, then $\lamax^\L(G)<\lamax^\Q(G)=\rho^\L(G)$.
If $\lamax^\L(G)<\rho^\L(G)$, it may occur  $\rho^\L(G) = \rho^\Q(G)$, which implies that the spectral radius is attained for some eigenvalue whose eigenvectors can not be scaled into $H$-eigenvectors, which are called {\it $N$-eigenvectors} of $\L(G)$.

In this paper we will discuss the above problem for the non-odd-bipartite generalized power hypergraphs $G^{k,{k \over 2} }$ constructed from non-bipartite simple graphs $G$, which will be introduced later.
In Section 2, we first give a method to compute the spectrum and the $H$-spectrum of $\L(G^{k,{k \over 2}})$ by computing the spectrum of certain matrices associated with
 the modified induced subgraph of the simple graph $G$.
 In particular, we given two explicit formulas for $\lamax^\L(G^{k,\frac{k}{2}})$ and $\rho^\L(G^{k,\frac{k}{2}})$ respectively.
 By using those results, in Section 3 we give a characterization for the equality $\rho^\L(G^{k,{k \over 2} }) = \rho^\Q(G^{k,{k \over 2} })$, i.e. $k$ is a multiple of $4$;in this case $\lamax^\L(G^{k,\frac{k}{2}})<\rho^\L(G^{k,\frac{k}{2}})$.
If $k\equiv 2 (\!\!\!\mod 4)$, then for sufficiently large $k$, $\lamax^\L(G^{k,\frac{k}{2}})<\rho^\L(G^{k,\frac{k}{2}})$.
 So, given a connected non-bipartite simple graph $G$, except a small number of $k$, we always have $\lamax^\L(G^{k,\frac{k}{2}})<\rho^\L(G^{k,\frac{k}{2}})$.
Motivated by the study of hypergraphs $G^{k,\frac{k}{2}}$, for a connected non-odd-bipartite hypergraph $G$,
we show that $Spec(\L(G))=Spec(\Q(G))$ (respectively, $Spec(\A(G))=-Spec(\A(G))$) if and only if $\L(G)$ and $\Q(G)$ (respectively, $\A(G)$ and $-\A(G)$) are similar via a complex (necessarily non-real) diagonal matrix with modular-$1$ diagonal entries.

In the paper \cite{SSW}, Shao et al. remarked that ``if $G$ is connected, then
 $$Hspec(\L(G))=Hspec(\Q(G)) \Longrightarrow Spec(\L(G))=Spec(\Q(G)).\eqno(1.2)$$
But we do not know whether the reverse implication is true or not.''
By our result, $Spec(\L(G))=Spec(\Q(G))$ is equivalent to that
$\L(G)$ is similar to $\Q(G)$ via a complex diagonal matrix with modular-$1$ diagonal entries.
However, by the results in \cite{SSW},  that $Hspec(\L(G))=Hspec(\Q(G))$ is equivalent to that
$\L(G)$ is similar to $\Q(G)$  via a diagonal matrix with $\pm 1$ diagonal entries.
So, if the complex diagonal matrix can be taken as real, then $Spec(\L(G))=Spec(\Q(G)) \Rightarrow  Hspec(\L(G))=Hspec(\Q(G))$.
But this happens only when $G$ is odd-bipartite by Theorem \ref{shao}.
Similar discussion can apply to $Spec(\A(G))$ and $Hspec(\A(G))$ for the spectral symmetric property.
So, for a connected non-odd-bipartite hypergraph $G$, the reverse implication in (1.2) is not true.

%

Finally we introduce the generalized power hypergraphs defined in \cite{KF}.

\begin{defi} {\em \cite{KF}}
Let $G=(V,E)$ be a simple graph. For any $k \ge 3$ and $1 \le s \le k/2$, the generalized power of $G$, denoted by $G^{k,s}$, is defined as
the $k$-uniform hypergraph with the vertex set $\{\v: v \in V\} \cup \{\mathbf{e}: e \in E\}$, and the edge set
$\{\u \cup \v \cup \mathbf{e}: e=\{u,v\} \in E\}$, where $\v$ is an $s$-set containing $v$ and $\mathbf{e}$ is a $(k-2s)$-set corresponding to $e$.
\end{defi}

Note that if $1 \le s<k/2$, then $G^{k,s}$ is a cored hypergraphs and hence is odd-bipartite. In particular, $G^{k,1}$ is exactly the $k$-th power of $G$.
If $s=k/2$ ($k$ being even), then $G^{k,s}$  is obtained from $G$ by only blowing up its vertices, $G^{2,1}=G$.
In this case, $\{u,v\}$ is an edge of $G$ if and only if $\u\cup\v$ is an edge of $G^{k,{k \over 2} }$, where we use the bold $\v$ to denote the blowing-up of the vertex $v$ in $G$.
For simplicity, we write $\u\v$ rather than $\u\cup\v$, and call $\u$ a {\it half edge} of $G^{k,{k \over 2} }$.

If $G=\Go$, a simple graph with loops (i.e. edges containing only one vertex), then $(\Go)^{k,s}$ will have loops containing $k-s$ vertices.
In particular, $(\Go)^{k,{k \over 2}}$ will have loops containing ${k \over 2}$ vertices.
That is, if $\{u\}$ is a loop of $\Go$, then the half edge $\u$ is a loop of $(\Go)^{k,{k \over 2}}$; see Fig. 1.1.

\begin{lemma}\label{NOB} {\em \cite{KF}} Let $G$ be a simple graph.
The hypergraph $G^{k,{k \over 2}}$ is non-odd-bipartite if and only if $G$ is non-bipartite.
\end{lemma}

\begin{lemma}\label{NOB2} {\em \cite{KF}} Let $G$ be a connected simple graph.
Then $\rho^\A(G)=\rho^\A(G^{k,{k \over 2}})$ and $\rho^\Q(G)=\rho^\Q(G^{k,{k \over 2}})$.
\end{lemma}

In the following for a simple graph $G$ and its generalized power hypergraph $G^{k,{k \over 2}}$,
 each vertex $u$ of $G$ is corresponding to the half edge $\u$ of $G^{k,{k \over 2}}$, and $u$ is always assumed to be contained in $\u$.
 Clearly, each vertex in $\u$ can be considered as $u$.
 In addition, all $k$-uniform hypergraphs are even uniform, i.e. $k$ is even.

\begin{center}
  \setlength{\unitlength}{1bp}%
  \begin{picture}(377.35, 306.95)(0,0)
  \put(0,0){\includegraphics{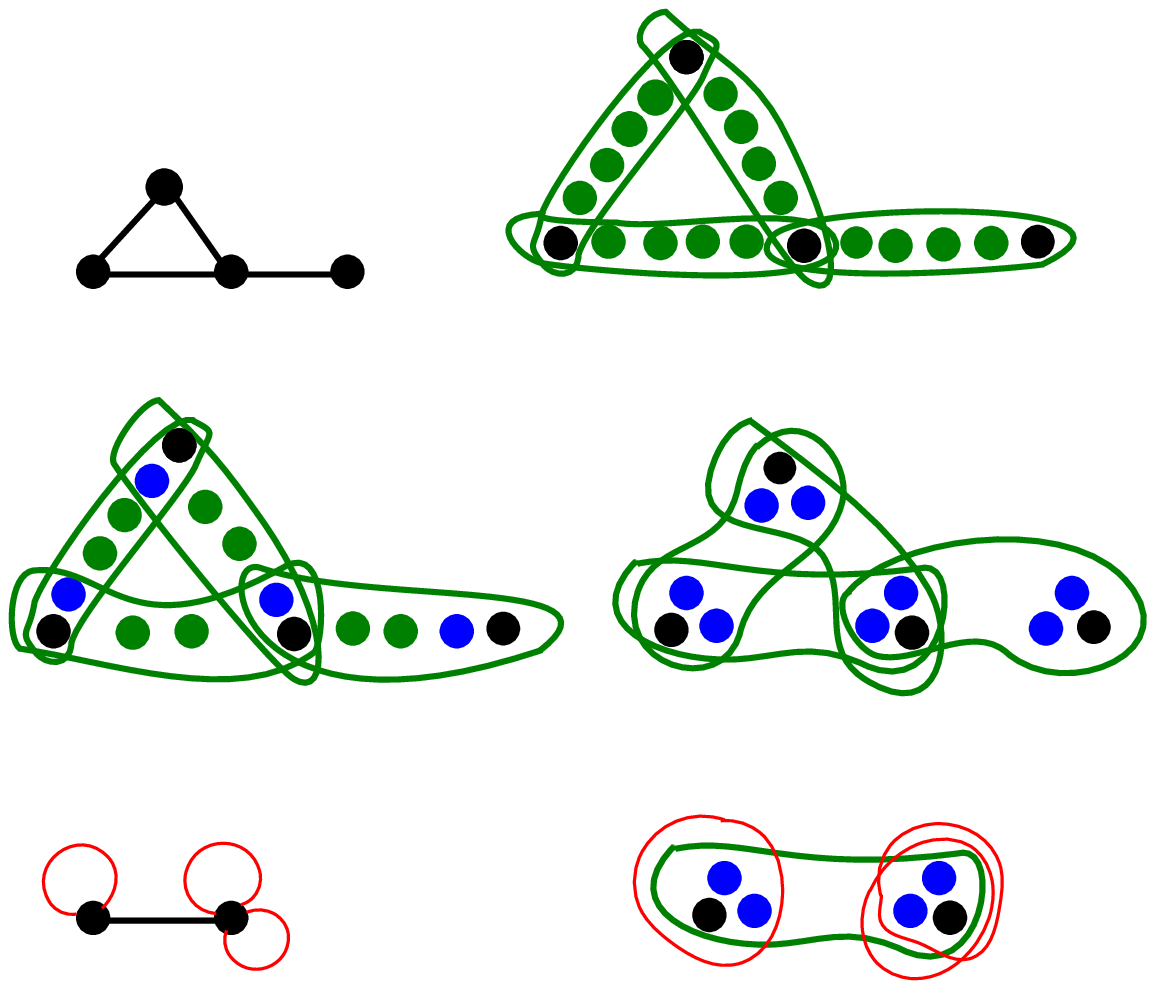}}
  \put(33.62,206.09){\fontsize{8.54}{10.24}\selectfont A simple graph $G$}
  \put(176.93,206.09){\fontsize{8.54}{10.24}\selectfont The power hypgergraph $G^6$}
  \put(5.67,90.43){\fontsize{8.54}{10.24}\selectfont The generalized power hypergraph $G^{6,2}$}
  \put(183.62,90.43){\fontsize{8.54}{10.24}\selectfont The generalized power hypergraph $G^{6,3}$}
  \put(26.49,234.10){\fontsize{8.54}{10.24}\selectfont $u$}
  \put(70.82,234.65){\fontsize{8.54}{10.24}\selectfont $w$}
  \put(26.49,48.00){\fontsize{8.54}{10.24}\selectfont $u$}
  \put(70.82,48.54){\fontsize{8.54}{10.24}\selectfont $w$}
  \put(5.67,7.51){\fontsize{8.54}{10.24}\selectfont Modified induced subgraph $G^{\large \mathtt  o}[u,w]$}
  \put(213.67,7.51){\fontsize{8.54}{10.24}\selectfont $G^{\large \mathtt  o}[u,w]^{k,\frac{k}{2}}$}
  \end{picture}%

\small{Fig. 1.1 (c.f. \cite{KF}) Constructing power hypergraphs from a simple graph, where a closed green curve represents an edge and a closed red curve represents a loop}
\end{center}

\section{The spectrum of Laplacian tensor}
In this section we will give a method to compute the spectra and the $H$-spectra of generalized power hypergraphs  $G^{k,{k \over 2}}$.
The eigenvector equation $\mathcal{L}(G)x^{k-1}=\lambda x^{[k-1]}$ could be interpreted as
$$ [d_v-\lambda] x_v^{k-1}= \sum_{\{v,v_2,v_3,\ldots, v_k\} \in E(G)} x_{v_2}x_{v_3} \cdots x_{v_k}, \mbox{~for each~} v \in V(G).\eqno(2.1)$$


Let $G$ be a simple graph on $n$ vertices possibly with loops.
Let $u$ be an arbitrary fixed vertex of $G^{k,\frac{k}{2}}$.
Define a vector $\x$ on $G^{k,\frac{k}{2}}$ such that $\x_u=1$ and $\x_v=0$ for any other vertices $v \ne u$.
It is easy to verify by (2.1) that $d_u$ is an eigenvalue (also an $H$-eigenvalue) of $\L(G^{k,\frac{k}{2}})$.

From the above fact, we find that the vertices in the same half edge of $G^{k,\frac{k}{2}}$ may have different values given by eigenvectors of $\L(G^{k,\frac{k}{2}})$.
However, if $\la \ne d_v$ for some vertex $v$,
we will have the following property on the eigenvectors associated with $\la$.

Denote by $d_\u$ the common degree of the vertices in $\u$.
For a nonempty subset $S \subseteq V(G^{k,{k \over 2}})$, denote $x^S:=\prod_{v \in S} x_v$, where $x$ is a vector defined on the vertices of $G^{k,{k \over 2} }$.

\begin{lemma}
\label{equal}Let $G$ be a simple graph possibly with loops.
Let $u$ and $\bar{u}$ be two vertices in the same half edge $\mathbf{u}$ of $G^{k,\frac{k}{2}}$.
If $\x$ is an eigenvector of $\L(G^{k,\frac{k}{2}})$ corresponding an eigenvalue $\lambda \ne d_\u$, then $\x_u^k=\x_{\bar{u}}^k$.
\end{lemma}

{\bf Proof:}
By the eigenvector equation (2.1),
$$ (d_\u-\lambda) \x_u^{k-1}=\sum_{\u\v \in E(G^{k,\frac{k}{2}})}\x^{\u \backslash \{u\}}\x^\v, ~~
(d_\u-\lambda) \x_{\bar{u}}^{k-1}=\sum_{\u\v \in E(G^{k,\frac{k}{2}})}\x^{\u \backslash \{\bar{u}\}}\x^\v.$$
So we have $(d_\u-\lambda) \x_u^{k}=(d_\u-\lambda) \x_{\bar{u}}^{k}$.
The result follows as $\lambda \ne d_\u$. \hfill $\blacksquare$

Let $\la$ be an eigenvalue of $\L(G^{k,\frac{k}{2}})$ such that $\la \notin \{d_\u: u \in V(G)\}$.
By Lemma \ref{equal}, the eigenvectors $\x$ of $\la$ have the common modulus on the vertices in each half edge $\u$, which will be denoted by $|\x_\u|$.
By Lemma \ref{equal}, if $\x_u=0$, then $\x_v=0$ for each $v \in \u$.
Otherwise, for each $v \in \u$,
$$\frac{\x_v}{\x_u}= e^{\i\frac{ 2 \pi \ell_{vu} }{k}}=:\E_{vu},\eqno(2.2)$$ where $\ell_{uu}=0$ and $\ell_{vu} \in \{0,1,\ldots,k-1\}$.
Suppose that $\x$ contains no zero entries. Define
$$\E_u:=\prod_{v \in \u} \E_{vu} = \frac{\x^{\u}}{\x_u^{k/2}}=e^{\frac{\i 2 \pi \sum_{v \in \u} \ell_{vu} }{k}}.\eqno(2.3)$$
Define $\mathcal{E}=\diag\{\E_u: u \in V(G)\}$, and
$$\L^\E(G)=\mathcal{D}(G)-\mathcal{E}\A(G)\mathcal{E}.\eqno(2.4)$$

If taking another vertex, say $\bar{u}$ as $u$, then $\E_{\bar{u}}=\pm \E_u$ as $\x^{\u}=\x_u^{k/2} \E_u=\x_{\bar{u}}^{k/2} \E_{\bar{u}}$ and $\x_u^k=\x_{\bar{u}}^k$.
Let $\mathcal{\bar{E}}=\diag\{\E_{\bar{u}}: u \in V(G)\}$. Then $\bar{\E}=\mathcal{E}S$, where $S$ is a diagonal matrix with $\pm 1$ on its diagonal.
So $\L^{\bar{\E}}(G)=\mathcal{D}(G)-\mathcal{\bar{E}}\A(G)\mathcal{\bar{E}}=S^{-1}\L^\E(G)S$, and hence $\L^{\bar{\E}}(G)$ has the same spectrum as $\L^\E(G)$.

\begin{lemma}\label{nonzero}
Let $\la$ be an eigenvalue of $\L(G^{k,\frac{k}{2}})$ corresponding to an eigenvector $\x$, where $G$ is a simple graph possibly with loops.
Suppose that $\x$ contains no zero entries.
If $\la \notin \{d_\u: u \in V(G)\}$ as an eigenvalue of $\L(G^{k,\frac{k}{2}})$,
then $\la$ is an eigenvalue of $\L^\E(G)$ with an eigenvector $x$ such that $x_u=\x_u^{k/2}$ for each $u \in V(G)$.
\end{lemma}

{\bf Proof:}
For each vertex $u \in \u$, $ (d_\u-\la)\x_u^{k-1}=\sum_{\u\w \in E(G^{k,\frac{k}{2}})} \x^{\u \backslash \{u\}}\x^\w.$
So by (2.3)
$$(d_\u-\la)\x_u^{k}=\sum_{\u\w \in E(G^{k,\frac{k}{2}})} \x^\u\x^\w=\sum_{\u\w \in E(G^{k,\frac{k}{2}})} \E_u \x_u^{k/2} \E_w \x_w^{k/2}.$$
$$(d_u-\la)\x_u^{k/2}=\sum_{uw \in E(G)} \E_u \x_w^{k/2} \E_w.$$
Therefore,
$\la$ is an eigenvalue of the matrix $\L^\E(G)$ with the eigenvector $x$ defined in the lemma. \hfill $\blacksquare$


The {\it modified induced subgraph} of a simple graph $G$ induced by the vertex subset $U \subseteq V(G)$, denoted by $\Go[U]$, is
the induced subgraph $G[U]$ together with $d_v(G)-d_v(G[U])$ loops on each vertex $v \in U$; see Fig. 1.1.
The Laplacian matrix of $\Go[U]$ is exactly $\L(G)[U]$, i.e. $\L(\Go[U])=\L(G)[U]$, and $\L(\Go[U]^{k,\frac{k}{2}})=\L(G^{k,\frac{k}{2}})[\mathbf{U}]$, where $\mathbf{U}=\cup\{\u: u\in U\}$.
Similarly, $\Q(\Go[U])=\Q(G)[U]$ and $\Q(\Go[U]^{k,\frac{k}{2}})=\Q(G^{k,\frac{k}{2}})[\mathbf{U}]$.

\begin{theorem} \label{main1}
Let $\la$ be an eigenvalue of $\L(G^{k,\frac{k}{2}})$ corresponding to an eigenvector $\x$,
where $G$ is a simple graph possibly with loops.
Suppose that $\la \notin \{d_\u: u \in V(G)\}$ as an eigenvalue of $\L(G^{k,\frac{k}{2}})$.
Let $\mathbf{U}=\cup\{\u: |\x_\u| > 0\}$ and $U=\{u: \u \subseteq \mathbf{U}\}$.
Let $\E=\diag\{\E_u: u \in U\}$ be defined as in (2.3).
Then the following results hold.

{\em (1)} $G^{k,\frac{k}{2}}[\mathbf{U}]$ contains no isolated half edges, and hence $G[U]$ contains no isolated vertices.

{\em (2)} $\la$ is an eigenvalue of $\L(G^{k,\frac{k}{2}})[\mathbf{U}]$ with $\x[\mathbf{U}]$ as an eigenvector.

{\em (3)} $\la$ is an eigenvalue of $\L^\E(\Go[U])$ with an eigenvector $x$ such that $x_u=\x_u^{k/2}$ for $u \in U$.

\end{theorem}

{\bf Proof:}
By (2.1), it is easy to verify the assertions (1) and (2).
Note that $\L(G^{k,\frac{k}{2}})[\mathbf{U}]=\L(\Go[U]^{k,\frac{k}{2}})$, the assertion (3) follows from Lemma \ref{nonzero} as $\x[\mathbf{U}]$ contains no zero entries.
\hfill $\blacksquare$

\begin{coro}\label{gou}
Each eigenvalue $\la$ of $\L(G^{k,\frac{k}{2}})$ is an eigenvalue of  $\L^\E(\Go[U])$ for some connected modified induced subgraph $\Go[U]$ and some choice of $\mathcal{E}$.
Furthermore, if $\la$ is an $H$-eigenvalue of $\L(G^{k,\frac{k}{2}})$, then $\la$ is an eigenvalue of $L(\Go[U])$.
\end{coro}

{\bf Proof:} Let $\x$ be an eigenvector of $\L(G^{k,\frac{k}{2}})$ corresponding to $\la$.
If $\la=d_\u$ for some half edge $\u$, then $\la$ is an $H$-eigenvalue of the Laplacian matrix $L(\Go[u])$.
Otherwise, let $\mathbf{U}$ and $U$ be defined as in Theorem \ref{main1}.
Then $\la$ is an eigenvalue of $\L^\e(\Go[U])$, where $\mathcal{E}=\diag\{\E_u: u \in U\}$.
We may assume that $\Go[U]$ is connected, as otherwise $\la$ must be a Laplacian eigenvalue of some connected component of $\Go[U]$.
If $\x$ is real, so is $\x[\mathbf{U}]$.
From the notations (2.2) and (2.3), for each half edge $\u \in \mathbf{U}$ and each vertex $v \in \u$, by Lemma \ref{equal}, $\E_{vu}=\pm 1$ and hence $\E_u=\pm 1$.
So, $\E=\E^{-1}$, and $$\L^\E(\Go[U])=\D(G)[U]-\E\A(G[U])\E=\E^{-1}(\D(G)[U]-\A(G[U]))\E=\E^{-1}L(\Go[U])\E,$$ which implies that
$\L^\E(\Go[U])$ has the same spectrum as $L(\Go[U])$.
Therefore $\la$ is an eigenvalue of $L(\Go[U])$.
\hfill $\blacksquare$

\begin{lemma}\label{contain} Let $G$ be a simple graph possibly with loops.
Let $\mathcal{E}=\diag \{\E_u: u \in V(G)\}$, where $\E_u=e^{\i\frac{ 2 \pi \ell_u }{k}}$ for some $\ell_u \in \{0,1,\ldots,k-1\}$.
Then each eigenvalue of $\L^\E(G)$ is an eigenvalue of $\L(G^{k,\frac{k}{2}})$.
\end{lemma}

{\bf Proof:}
Let $\la$ be an eigenvalue of $\L^\e(G)$ associated with the eigenvector $x$.
For each half edge $\u$ of $G^{k,\frac{k}{2}}$,
there exists a function $f_\u: \u \to \{0,1,\ldots,k-1\}$ such that $f_\u(u)=0$ and $e^{\i\frac{ 2 \pi \sum_{v \in \u} f_\u(v) }{k}}=\E_u$.
Now define a vector $\x$ defined on $G^{k,\frac{k}{2}}$ such that for each half edge $\u$ and each $v \in \u$,
$$\x_v=x_u^{2/k}e^{\i\frac{ 2 \pi f_\u(v) }{k}},\eqno(2.5)$$
where $x_u^{2/k}$ is a root of the equation $\alpha^{k/2}=x_u$.
By the eigenvector equation of $\L^\E(G)$, for each vertex $u$,
$$ (d_u -\la)x_u=\sum_{uw \in E(G)} \E_u x_w \E_w.$$
So
$$ (d_u -\la)\x_u^{k/2}=\sum_{uw \in E(G)} \E_u \x_w^{k/2} \E_w.$$
$$ (d_\u -\la)\x_u^{k-1}=\sum_{uw \in E(G)} \E_u \x_u^{k/2-1} \x_w^{k/2} \E_w=\sum_{\u\w \in E(G^{k,\frac{k}{2}})} \x^{\u \backslash \{u\}} \x^\w.$$
For any other vertex $v \in \u$,
$$ (d_\u -\la)\x_v^{k-1}=(d_\u -\la)(\x_u e^{\frac{\i 2 \pi f_\u(v) }{k}})^{k-1} =\sum_{\u\w \in E(G^{k,\frac{k}{2}})} \x^{\u \backslash \{u\}} \x^\w e^{-\frac{\i 2 \pi f_\u(v) }{k}}=\sum_{\u\w \in E(G^{k,\frac{k}{2}})} \x^{\u \backslash \{v\}} \x^\w.$$
Therefore $\la$ is an eigenvalue of $\L(G^{k,\frac{k}{2}})$ with the eigenvector $\x$ defined as in (2.5).
\hfill $\blacksquare$


By Lemma \ref{contain}, if taking $\mathcal{E}=\mathcal{I}$, then each eigenvalue of $\L(G)$ is an eigenvalue of $\L(G^{k,\frac{k}{2}})$.
We will show those eigenvalues of $\L(G)$ are really $H$-eigenvalue of $\L(G^{k,\frac{k}{2}})$.

\begin{lemma}\label{contain-H} Let $G$ be a simple graph possibly with loops.
Each eigenvalue of $\L(G)$ is an $H$-eigenvalue of $\L(G^{k,\frac{k}{2}})$.
\end{lemma}

{\bf Proof:}
Let $x$ be an eigenvector of $\L(G)$ corresponding to an eigenvalue $\lambda$.
Let $\x$ be a vector defined on $G^{k,\frac{k}{2}}$ as follows. For each $ u \in V(G)$,
$$ \x_u=\sgn(x_u)|x_u|^{2/k}, ~~\x_v=|x_\u|^{2/k}, \mbox{~for each vertex~} v \in \u \backslash \{u\}.\eqno(2.6)$$
Then $$\x^{\u}=x_u, \mbox{~for each~} u \in V(G).$$
Also, since $k$ is even,
$$\x_u^{k-1}=(\sgn(x_u)|x_u|^{2/k})^{k-1}=\sgn(x_u)|x_u|(|x_u|^{2/k})^{k/2-1}=x_u \x^{\u\backslash \{u\}}.$$
By the eigenvector equation of $\L(G)$,
$(d_u-\lambda) x_u=\sum_{uw \in E(G)}x_w$, so we have
$$(d_u-\lambda) \x_u^{k-1} =(d_u-\lambda)x_u \x^{\u\backslash \{u\}}=\sum_{uw \in E(G)}x_w \x^{\u\backslash \{u\}}=\sum_{\u\w \in E(G^{k,\frac{k}{2}})}\x^{\u \backslash \{u\}}\x^\w.$$
For any other vertex $v \in \u$,
\begin{align*}
(d_u-\lambda) \x_v^{k-1} &=(d_u-\lambda)(\sgn(x_u) \x_u)^{k-1}=\sgn(x_u) (d_u-\lambda)\x_u^{k-1}\\
&=\sum_{\u\w \in E(G^{k,\frac{k}{2}})}\sgn(x_u) \x^{\u \backslash \{u\}}\x^\w=
\sum_{\u\w \in E(G^{k,\frac{k}{2}})}\x^{\u \backslash \{v\}}\x^\w.
\end{align*}
So $\lambda$ is an $H$-eigenvalue of $\L(G^{k,\frac{k}{2}})$.
\hfill $\blacksquare$

\begin{coro}\label{Gou}
Let $G$ be a simple graph, and let $\Go[U]$ be a connected modified induced subgraphs of $G$.
Let $\mathcal{E}=\diag \{\E_u: u \in U\}$, where $\E_u=e^{\i\frac{ 2 \pi \ell_u }{k}}$ for some $\ell_u \in \{0,1,\ldots,k-1\}$.
Then each eigenvalue of $\L^\E(\Go[U])$ is an eigenvalue of $\L(G^{k,\frac{k}{2}})$.
In particular, each eigenvalue of $\L(\Go[U])$ is an $H$-eigenvalue of $\L(G^{k,\frac{k}{2}})$.
\end{coro}

{\bf Proof:}
By Lemma \ref{contain}, if $\la$ is an eigenvalue of $\L^\E(\Go[U])$ with an eigenvector $x$, then $\la$ is an eigenvalue of
   $\L(\Go[U]^{k,\frac{k}{2}})=\L(G^{k,\frac{k}{2}})[\mathbf{U}]$ with an eigenvector $\x$ whose entries are defined as in (2.5), where $\mathbf{U}=\cup\{\u: u \in U\}$.
Extending the eigenvector $\x$ defined on $\Go[U]^{k,\frac{k}{2}}$ to $G^{k,\frac{k}{2}}$ by assigning zeros to the vertices outside $\mathbf{U}$,
we will get a vector $\mathbf{y}$.
It is easy to verify by (2.1) that $\mathbf{y}$ is an eigenvector of $\L(G^{k,\frac{k}{2}})$ corresponding the eigenvalue $\la$.

If $\E=\mathcal{I}$, $\L^\E(\Go[U])=\L(\Go[U])$ and $x$ could be taken real.
In this case, by Lemma \ref{contain-H} we take the real eigenvector $\x$ whose entries are defined as in (2.6).
Then by a similar discussion, $\la$ is an $H$-eigenvalue of $\L(G^{k,\frac{k}{2}})$. \hfill $\blacksquare$

By Corollary \ref{gou} and Corollary \ref{Gou}, we get the following main result.

\begin{theorem}\label{main2}
Let $G$ be a simple graph.
Then, regardless of multiplicities, the spectrum of $\L(G^{k,\frac{k}{2}})$
 consists of all eigenvalues of $\L^\E(\Go[U])$ for all choices of $\mathcal{E}$ as defined in Corollary \ref{Gou} and all connected modified induced subgraphs $\Go[U]$ of $G$.

 Furthermore, regardless of multiplicities, the $H$-spectrum of $\L(G^{k,\frac{k}{2}})$ consists of all eigenvalues of $\L(\Go[U])$ for all connected modified induced subgraphs $\Go[U]$ of $G$.
\end{theorem}

\begin{coro}\label{LE relation}
Let $G$ be a simple graph.
Then $\lamax^\L(G^{k,\frac{k}{2}})=\lamax^\L(G)$, $\rho^\L(G^{k,\frac{k}{2}})=\max\{\rho(\L^\E(\Go[U]))\}$,
where the maximum is taken over all all choices of $\mathcal{E}$ as defined in Corollary \ref{Gou} and all connected modified induced subgraphs $\Go[U]$ of $G$.
\end{coro}

{\bf Proof:}
By the interlacing theorem of the eigenvalues of real symmetric matrices (see \cite{hc}), $\lamax^\L(G)$ is the maximum of all largest eigenvalues of the principal submatrices of $\L(G)$. The first equality follows from Theorem \ref{main2}.
The second equality is easily seen also by Theorem \ref{main2}.
\hfill $\blacksquare$

Along the line of discussion in this section, one can easily get the spectrum of the adjacency tensor or the signless Laplacian tensor, where the $H$-spectra of these tensors
are discussed in \cite{KF2}.

\begin{theorem}
Let $G$ be a simple graph.
Then, regardless of multiplicities, the spectrum of $\A(G^{k,\frac{k}{2}})$ (respectively, $\Q(G^{k,\frac{k}{2}})$)
 consists of all eigenvalues of $\A^\E(G[U])$ (respectively, $\Q^\E(\Go[U])$) for all choices of $\mathcal{E}$ as defined in Corollary \ref{Gou} and all connected induced subgraphs $G[U]$ (respectively, all connected modified induced subgraphs $\Go[U]$) of $G$, where $\A^\E(G[U])=\E \A(G[U]) \E$ and $\Q^\E(\Go[U])=\D(G)[U]+\E \A(G[U]) \E$.

 Furthermore, regardless of multiplicities, the $H$-spectrum of $\A(G^{k,\frac{k}{2}})$ (respectively, $\Q(G^{k,\frac{k}{2}})$) consists of all eigenvalues of $\A(G[U])$ (respectively, $\Q(\Go[U])$) for all connected induced subgraphs $G[U]$ (respectively, all connected modified induced subgraphs $\Go[U]$) of $G$.
\end{theorem}

\section{The largest $H$-eigenvalue and spectral radius of Laplacian tensor}
Let $G$ be a connected simple graph.
If $G$ is bipartite, then $G^{k,\frac{k}{2}}$ is odd-bipartite by Lemma \ref{NOB}.
So by Theorem \ref{Hu}, $\lamax^\L(G^{k,\frac{k}{2}})=\lamax^\Q(G^{k,\frac{k}{2}})$, which implies that
$$\rho^\L(G^{k,\frac{k}{2}})=\lamax^\L(G^{k,\frac{k}{2}})=\lamax^\Q(G^{k,\frac{k}{2}})=\rho^\Q(G^{k,\frac{k}{2}}).$$
If $G$ is non-bipartite, then $G^{k,\frac{k}{2}}$ is non-odd-bipartite also by Lemma \ref{NOB}.
By Theorem \ref{Hu}, $$\lamax^\L(G^{k,\frac{k}{2}})<\lamax^\Q(G^{k,\frac{k}{2}})=\rho^\Q(G^{k,\frac{k}{2}}).$$
However, it may occur that $\rho^\L(G^{k,\frac{k}{2}})=\rho^\Q(G^{k,\frac{k}{2}})$.

\begin{lemma}\label{equalspec}
Let $G$ be a connected non-bipartite graph.
If $k$ is a multiple of $4$, then $\rho^\L(G^{k,\frac{k}{2}})=\rho^\Q(G^{k,\frac{k}{2}})$, or equivalently $Spec(\L(G^{k,\frac{k}{2}}))=Spec(\Q(G^{k,\frac{k}{2}}))$.
\end{lemma}

{\bf Proof:}
It suffices to prove that $\rho^\Q(G^{k,\frac{k}{2}})$ is an eigenvalue of $\L(G^{k,\frac{k}{2}})$ as  $\rho^\L(G^{k,\frac{k}{2}}) \le \rho^\Q(G^{k,\frac{k}{2}})$.
Let $\x$ be an eigenvector $\Q(G^{k,\frac{k}{2}})$ corresponding to $\rho^\Q(G^{k,\frac{k}{2}})=:\rho$.
By the eigenvector equation of $\Q(G^{k,\frac{k}{2}})$, for each vertex $u \in \u$,
$$(\rho-d_\u)\x_u^{k-1}=\sum_{\u\w \in E(G^{k,\frac{k}{2}})}\x^{\u\backslash \{u\}}\x^\w.\eqno(3.1)$$
Define a vector $\y$ such that for each half edge $\u$,
$$\y_u=\i \x_u,  \y_v=\x_v \mbox{~ for any other~} v \in \u \backslash \{u\}.\eqno(3.2)$$
Noting that $k$ is a multiple of $4$, by (3.1) it is easy to verify
$$(d_\u-\rho)\y_u^{k-1}=\sum_{\u\w \in E(G^{k,\frac{k}{2}})}\y^{\u\backslash \{u\}}\y^\w,$$
and for any other vertex $v \in \u$,
$$(d_\u-\rho)\y_v^{k-1}=\sum_{\u\w \in E(G^{k,\frac{k}{2}})}\y^{\u\backslash \{v\}}\y^\w.$$
So $\rho^\Q(G^{k,\frac{k}{2}})$ is an eigenvalue of $\L(G^{k,\frac{k}{2}})$ with $\y$ as an eigenvector. \hfill $\blacksquare$

We give some remarks for Lemma \ref{equalspec}.
For each half edge $\u$ of $G^{k,\frac{k}{2}}$, define
$$\Gamma_u=\i, \Gamma_v=1, \mbox{~for any other vertex~} v \in \u.$$
Then we get a diagonal matrix $\Gamma=\diag\{\Gamma_v:  v \in V(G^{k,\frac{k}{2}})\}$.
From the proof of Lemma \ref{equalspec}, if $\x$ is an eigenvector of $\Q(G^{k,\frac{k}{2}})$ corresponding to an eigenvalue $\la$,
then $\Gamma \x$ is an eigenvector of $\L(G^{k,\frac{k}{2}})$ also corresponding to the eigenvalue $\la$.
Furthermore, according the tensor product introduced in \cite{Shao}, $$\L(G^{k,\frac{k}{2}})=\Gamma^{-(k-1)}\Q(G^{k,\frac{k}{2}})\Gamma,\eqno(3.3)$$
which implies that $\L(G^{k,\frac{k}{2}})$ and $\Q(G^{k,\frac{k}{2}})$ are diagonal similar, and hence
$Spec(\L(G^{k,\frac{k}{2}}))=Spec(\Q(G^{k,\frac{k}{2}}))$ by \cite[Theorem 2.3]{Shao} though $Hspec(\L(G^{k,\frac{k}{2}}))\ne Hspec(\Q(G^{k,\frac{k}{2}}))$
by Theorem \ref{shao} as $G^{k,\frac{k}{2}}$ is not odd-bipartite.
From (3.3), one can get
$$-\A(G^{k,\frac{k}{2}})=\Gamma^{-(k-1)}\A(G^{k,\frac{k}{2}})\Gamma,\eqno(3.4)$$
so $Spec(\A(G^{k,\frac{k}{2}}))=-Spec(\A(G^{k,\frac{k}{2}}))$, i.e. the spectrum is symmetric with respect to the origin, though $Hspec(\A(G^{k,\frac{k}{2}}))=-Hspec(\A(G^{k,\frac{k}{2}}))$ by Theorem \ref{shao} as $G^{k,\frac{k}{2}}$ is not odd-bipartite.

Secondly, the eigenvector $\y$ in the proof of Lemma \ref{equalspec} can also be defined in a way different from (3.2).
For each half edge $\u$, arbitrarily choose $\frac{k}{4}$-subset $U$ from $\u$,
and define $\y_u=e^{\i \frac{2\pi}{k}}$ if $u \in U$, and $\y_v =\x_v$ if $v \in \u \backslash U$.
One can also find a diagonal matrix $\Gamma$ based on this definition of $\y$ to make (3.3) and (3.4) hold.

Motivated by the above discussion, we get a result complementary to Theorems \ref{shao0} and \ref{shao}.

\begin{theorem} \label{fan}
Let $G$ be a connected non-odd-bipartite even uniform hypergraph.
Then the following are equivalent.

{\em(1)} $\rho^\L(G)=\rho^\Q(G)$.

{\em (2)} $\L(G)$ and $\Q(G)$ are similar via a complex (necessarily non-real) diagonal matrix with modular-$1$ diagonal entries.

{\em(3)} $Spec(\L(G))=Spec(\Q(G))$.

{\em (4)} $\A(G)$ and $-\A(G)$ are similar via a complex (necessarily non-real) diagonal matrix with modular-$1$ diagonal entries.

{\em(5)} $Spec(\A(G))=-Spec(\A(G))$.

{\em(6)} $-\rho^\A(G) \in Spec(\A(G))$.

\end{theorem}

{\bf Proof:} It is clear that (2) $\Rightarrow$ (3) $\Rightarrow$ (1) and (4) $\Rightarrow$ (5) $\Rightarrow$ (6) by \cite[Theorem 2.3]{Shao}.
We will take the proof technique from \cite{SSW}.
If $\rho^\L(G)=\rho^\Q(G)$, taking $\la=\rho^\Q(G)e^{\i \phi}$ as an eigenvalue of $\L(G)$, by Perron-Frobenius Theorem for nonnegative weakly irreducible tensors (see \cite{YY3}), there exists a nonsingular diagonal matrix $\Gamma$ with $|\Gamma|=\mathcal{I}$  such that
$$\L(G)=e^{\i \phi} \Gamma^{-(k-1)}\Q(G)\Gamma.\eqno(3.5)$$
So, $e^{\i \phi}=1$ by comparing the diagonal entries of both sides of (3.5), and
$$\L(G)=\Gamma^{-(k-1)}\Q(G)\Gamma,\eqno(3.6)$$
From (3.6) we have
 $$-\A(G)=\Gamma^{-(k-1)}\A(G)\Gamma.$$
  So, if (1) holds, we can get (2) and (4).
 Note that the matrix  $\Gamma$ can not be taken as real; otherwise, $\Gamma$ would have both $1$ and $-1$ along its diagonal, and then $G$ is odd-bipartite by \cite[Theorem 2.1]{SSW}; a contradiction.

Now suppose (6) holds, i.e. $-\rho^\A(G) \in Spec(\A(G))$. By Perron-Frobenius Theorem, there also exists a nonsingular diagonal matrix $\bar{\Gamma}$ with $|\bar{\Gamma}|=\mathcal{I}$ such that
$$\A(G)=-\bar{\Gamma}^{-(k-1)}\A(G)\bar{\Gamma},\eqno(3.7)$$
where the matrix $\bar{\Gamma}$ can not be taken as real by a similar discussion as the above.
From (3.7) we have
 $$\L(G)=\bar{\Gamma}^{-(k-1)}\Q(G)\bar{\Gamma},$$
 which implies that (2) holds. \hfill $\blacksquare$

From the proof of Theorem \ref{fan}, that $Spec(\L(G))=Spec(\Q(G))$ is equivalent to that
$\L(G)$ is similar to $\Q(G)$ via a complex diagonal matrix with modular-$1$ diagonal entries.
However, by the results in \cite{SSW},  that $Hspec(\L(G))=Hspec(\Q(G))$ is equivalent to that
$\L(G)$ is similar to $\Q(G)$  via a diagonal matrix with $\pm 1$ diagonal entries.
So, if the complex diagonal matrix can be taken as real, then $Spec(\L(G))=Spec(\Q(G)) \Rightarrow  Hspec(\L(G))=Hspec(\Q(G))$.
But this happens only when $G$ is odd-bipartite by Theorem \ref{shao}.
Similar discussion can apply to $Spec(\A(G))$ and $Hspec(\A(G))$ for the spectral symmetric property.

\begin{theorem} \label{main3}
Let $G$ be a connected non-bipartite graph.
Then $\rho^\L(G^{k,\frac{k}{2}})=\rho^\Q(G^{k,\frac{k}{2}})$ if and only if $k$ is a multiple of $4$.
In this case, $\lamax^\L(G^{k,\frac{k}{2}})<\rho^\L(G^{k,\frac{k}{2}})$.
\end{theorem}

{\bf Proof:}
The sufficiency follows by Lemma \ref{equalspec}.
By Corollary \ref{LE relation}, suppose that $\rho^\L(G^{k,\frac{k}{2}})=\rho(\L^\E(\Go[U]))$ for some connected modified induced subgraphs $\Go[U]$ of $G$ and some $\mathcal{E}$.
As $|\L^\E(\Go[U])|=\Q(\Go[U])$, by Perron-Frobenius Theorem for nonnegative weakly irreducible tensors (see \cite{YY3}) or for nonnegative irreducible matrices (see \cite{hc}) and Lemma \ref{NOB2},
$$ \rho^\L(G^{k,\frac{k}{2}})=\rho(\L^\E(\Go[U])) \le \rho(|\L^\E(\Go[U])|)= \rho^\Q(\Go[U]) \le \rho^Q(G)=\rho^\Q(G^{k,\frac{k}{2}}).$$
If $\rho^\L(G^{k,\frac{k}{2}})=\rho^\Q(G^{k,\frac{k}{2}})$, then $\rho^\Q(\Go[U]) = \rho^Q(G)$, which implies that $U=V(G)$ as $G$ is connected.
So $\rho(\L^\E(G))=\rho^\Q(G)$.
Assume that $\la=e^{\i \phi}\rho^\Q(G)$ is an eigenvalue of $\L^\E(G)$.
By Perron-Fronenius Theorem, there exists a diagonal matrix
$\Gamma=\diag \{ e^{\i \theta_u}: u \in V(G) \}$ such that $$ \L^\E(G)=e^{\i \phi} \Gamma^{-1}\Q(G) \Gamma. \eqno(3.8)$$
From (3.8) we have
$$e^{\i \phi}\Gamma^{-1}\D(G)\Gamma=\D(G), ~e^{\i \phi}\Gamma^{-1}\E\A(G)\E\Gamma=-\A(G). \eqno(3.9)$$
So, $e^{\i \phi}=1$.
As $G$ is non-bipartite, letting $C_{2m+1}$ be an odd cycle of $G$ with edges $v_iv_{i+1}$ for $i=1,2,\ldots,2m+1$, where $v_{2m+2}=v_1$.
Using the second equality of (3.9), for $i=1,2,\ldots,2m+1$,
$$e^{-\i \theta_{v_i}} \E_{v_i} \E_{v_{i+1}} e^{\i \theta_{v_{i+1}}}=-1.$$
Thus
$$ \prod_{i=1}^{2m+1} \left( e^{-\i \theta_{v_i}} \E_{v_i} \E_{v_{i+1}} e^{\i \theta_{v_{i+1}}}\right)=-1,$$
and hence $$\prod_{i=1}^{2m+1} \E_{v_i}^2=-1.$$
Noting that $\E_v=e^{\i\frac{ 2 \pi \ell_u }{k}}$ for some $\ell_u \in \{0,1,\ldots,k-1\}$,
$$ e^{\i\frac{4 \pi \sum_{i=1}^{2m+1} \ell_{v_i} }{k}}=-1,$$
which implies that $k$ is a multiple of $4$.\hfill $\blacksquare$

Next we discuss the case of $k \equiv 2 (\!\!\!\mod 4)$.
In this case, $\rho^\L(G^{k,\frac{k}{2}}) < \rho^\Q(G^{k,\frac{k}{2}})$ by Theorem \ref{main3}.
But, can we have $\lamax^\L(G^{k,\frac{k}{2}}) =\rho^\L(G^{k,\frac{k}{2}})$?

\begin{theorem}\label{main4}
Let $G$ be a connected non-bipartite graph. Suppose that $k\equiv 2 (\!\!\!\mod 4)$.
Then for sufficiently large $k$, $\lamax^\L(G^{k,\frac{k}{2}})<\rho^\L(G^{k,\frac{k}{2}})$.
\end{theorem}

{\bf Proof:}
Let $k=4l+2$, and let $\tilde{\E}=e^{\i \frac{2\pi l}{k}} \mathcal{I}$.
Then $$\L^{\tilde{\E}}(G)=\D(G)-\tilde{\E}\A(G)\tilde{\E}=\D(G)-e^{\i \frac{2\pi l}{2l+1}}\A(G).$$
If $k \to \infty$ (i.e. $l \to \infty$), then $ \L^{\tilde{\E}}(G) \to \D(G)+\A(G)=\Q(G)$.
As $\rho(\L^{\tilde{\E}}(G))$ is continuous in the entries of $\L^{\tilde{\E}}(G)$, if $k \to \infty$,
$$\rho(\L^{\tilde{\E}}(G)) \to \rho(\Q(G))=\rho^\Q(G^{k,\frac{k}{2}}). $$
By Corollary \ref{LE relation},
$$\rho^\L(G^{k,\frac{k}{2}})=\max\{\rho(\L^\E(\Go[U]))\} \ge \rho(\L^{\tilde{\E}}(G)).$$
Note that $\rho^\L(G^{k,\frac{k}{2}}) < \rho^\Q(G^{k,\frac{k}{2}})$ by Theorem \ref{main3}.
So, $$\rho^\L(G^{k,\frac{k}{2}}) \to \rho^\Q(G^{k,\frac{k}{2}})=\rho(\Q(G)).\eqno(3.10)$$
Since $G$ is non-bipartite, by Corollary \ref{LE relation}, $$\lamax^\L(G^{k,\frac{k}{2}})=\lamax^\L(G)=\rho(\L(G)) < \rho(\Q(G)).\eqno(3.10)$$
Combining (3.10) and (3.11), for sufficiently large $k$, $\lamax^\L(G^{k,\frac{k}{2}})<\rho^\L(G^{k,\frac{k}{2}})$.\hfill $\blacksquare$

By Theorem \ref{main3} and Theorem \ref{main4}, we pose the following conjecture.

\begin{conj} Let $G$ be a connected non-odd-bipartite hypergraph. Then $\lamax^\L(G)<\rho^\L(G)$.
\end{conj}

For a connected non-odd-bipartite hypergraph $G$, by Theorem \ref{Hu}, $\lamax^L(G)< \lamax^Q(G)=\rho^\Q(G)$.
If $\rho^L(G)=\rho^\Q(G)$, surely,  $\lamax^\L(G)<\rho^\L(G)$, and the above conjecture holds.
So, it suffices to consider those hypergraphs $G$ with $\rho^L(G)<\rho^\Q(G)$ for the conjecture.


\end{document}